\documentstyle[12pt]{article}

\begin{document}
\begin{flushright}
 q-alg/9705013 \\
\end{flushright}
\begin{center}
{\LARGE
BRST Algebra Quantum Double and Quantization of the Proper Time
Cotanjent Bundle
} \\[1cm]
{\bf Lyakhovsky V.D.\\
Theoretical Department \\
Institute of Physics \\
St.Petersburg State University \\
198904 St.Petersburg \\
Russia} \\
e-mail: lyakhovs@snoopy.phys.spbu.ru \\[10mm]
{\bf Tkach V.I.\\
Instituto de Fisica \\
Universidad de Guanajuato\\
Lomas del Bosque 103 \\
Col. Lomas del Campestre, A.P. E-143, C.P. 37150 \\
Leon, Gto. \\
Mexico} \\
e-mail: vladimir@ifug1.ugto.mx \\[15mm]

{\bf Abstract}

\end{center}

The quantum double for the quantized $BRST$ superalgebra is studied. The
corresponding $R$-matrix is explicitly constucted. The Hopf algebras of the
double form an analytical variety with coordinates
described by the canonical deformation parameters. This provides the
possibility to construct the nontrivial quantization of the proper time
supergroup cotangent bundle. The group-like classical limit for this
quantization corresponds to the generic super Lie bialgebra of the double.

\vspace{2cm}

\section{Introduction}
The theories of Casalbuoni-Brink-Schwarz (CBS) superparticle \cite{A} are
fundamentally related to supersymmetric field theories and strings.
Superparticle orbits are determined up to local fermionic (Siegel)
transformations \cite{B} , which play crutial role in removing the
unphysical degrees of freedom. For the case of superparticle it has
been shown \cite{C} that Siegel symmetry can be interpreted as the
usual local proper-time supersymmetry (PTSA). The equivalence
between CBS-superparticle and the spinning particle was established \cite{D}
by identifying Lorentz-covariant Siegel generator with the local proper-time
supersymmetry of the spinning particle \cite{E}.

To quantize such models it is natural to apply the BRST formalism, which is
manifestly Lorentz-invariant. For the point particle case the BRST
quantization starts with the Faddeev-Popov prescription and the extraction
of a new nilpotent symmetry operator. The latter can be included in the
algebra ILI(1) \cite{F}.

Thus the symmetry algebra of a system with superparticles contain
both $BRST$ and $PTSA$ subalgebras. The simpliest possible unification
of them is the direct sum. It is natural to consider the properties of
quantum analogues of $(PTSA)
\oplus (BRST)$.
On the other hand $BRST$ algebra itself can be treated as a deformation of
the trivial algebra of coordinate functions for the superparticle. So one
can equally consider $q$-deformations of a unification of $PTSA$ with
Abelian superalgebra creating the $BRST$ subalgebra in the process of
deformation. In this case the initial unification is a semidirect sum
corresponding to the coadjoint action.

The significant feature of the symmetries $PTSA$ and $BRST$ is that their
superalgebras are dual. This gives the opportunity to obtain the
necessary $q$-deformed symmetry by constructing Drinfeld double for a
quantized $(PTSA)_q$ superalgebra. The latter
is easily obtained using the method developed in \cite{Kulish}.

In this paper we demonstrate that the Hoph algebra of the quantum double
$SD(PTSA_q, BRST_q)$ can be treated as a quantized symmetry for both
interpretation schemes presented above. For the first one the double must
be considered as a quantum group corresponding to the algebra $(PTSA)\oplus
(BRST)^{\rm opp}$. In the second approach the multiplications in $SD$ are
treated as the deformedalgebraof coadjoint extension of $(PTSA)$.

The paper is organized as follows. In the second section all the necessary
algebraic constructions are obtained including the explicite expression of
${\cal R}$-matrix for $SD(PTAS_q,BRST_q)$. In section 3 the dual canonical
parameters are introduced in $SD$. This gives the possibility to construct
the limit transitions connecting different Poisson structures in the
created set of Hopf algebras. All the necessary classical limits are
explicitely realized. The obtained results are discussed in section 3 from
the point of view of possible physical interpretation.

\section{The BRST algebra quantum double}

Let the Hopf algebra with the generators $\left\{ T,S\right\} $ and the
defining relations
\begin{equation}
\label{ptsa-def}
\begin{array}{l}
[T,S]=0; \\
\left\{ S,S\right\} =2
\frac{\sinh \left( hT\right) }{\sinh \left( h\right) }; \\ \Delta T=T\otimes
1+1\otimes T; \\
\Delta S=e^{hT/2}\otimes S+S\otimes e^{-hT/2}.
\end{array}
\end{equation}
be interpreted as the proper-time quantum superalgebra $(PTSA_q)$. Chose the
following quantization of the two-dimensional $BRST$-algebra with basic
elements $\left\{ \tau ,\xi \right\} $:
\begin{equation}
\label{brst-def}
\begin{array}{l}
[\tau ,\xi ]=\frac h2\xi ; \\
\left\{ \xi ,\xi \right\} =0; \\
\Delta \tau =\tau \otimes 1+1\otimes \tau +\frac h{\sinh \left( h\right)
}\xi \otimes \xi ; \\
\Delta \xi =\xi \otimes 1+1\otimes \xi .
\end{array}
\end{equation}

Consider generators $\tau $ and $\xi $ as dual to $T$ and $S$.
Then the algebra (\ref{brst-def}) can be treated as dual
opposite to $PTSA_q$, that is the
$(PTSA_q)^*$ , with opposite comultiplication and inverse antipode.

Note that according to the quantum duality
principle \cite{Drin,Sem} the $PTSA_q$ algebra
defines also the quantization of the 2-dimensional vector quantum group
described by the coproducts in (\ref{ptsa-def}). This is the semidirect
product of two abelian groups and its supergroup nature is reflected only by
the fact that its topological space is a superspace. The quantum supergroup
(different from the previous one) is also defined by the Hopf algebra $BRST
_q $ (see $\Delta$'s in (\ref{brst-def})).

To obtain the quantum superdouble $SD(PTSA_q,BRST_q)$ one can start by
constructing the corresponding universal element. Let us define
the Poincare-Birkhoff-Witt-basis for $PTSA_q$ and $BRST_q$:
\begin{equation}
\label{basis}
\begin{array}{l}
1,\xi ,
\frac{\tau ^n}{n!},\frac{\xi \tau ^n}{n!} \\ 1,S,\frac{T^n}{n!},\frac{ST^n}{
n!}.
\end{array}
\end{equation}
The universal element can be written in the form
\begin{equation}
\label{rmat}{\cal R}=(1\otimes 1+S\otimes \xi )e^{T\otimes \tau }.
\end{equation}
Its main properties are easily checked with the help of an auxiliary relation
\begin{equation}
\label{supprel}
\begin{array}{c}
\left( 1\otimes 1\otimes 1+
\frac{(e^{2hT}-1)}{e^h-e^{-h}}\otimes \xi \otimes \xi \right) \exp \left(
T\otimes 1\otimes \tau +T\otimes \tau \otimes 1\right) = \\ =\exp \left(
T\otimes \tau \otimes 1+T\otimes 1\otimes \tau +\frac h{\sinh \left(
h\right) }T\otimes \xi \otimes \xi \right) .
\end{array}
\end{equation}

Next step should involve the construction of the multiplication rules
consistent with this ${\cal R}$-matrix. For any pair of
dual Hopf algebras $H$ and $H^{*}$ with the basic elements $\{e_s\}$ and $
\{e^t\}$ and the universal element ${\cal R}=e_s\otimes e^s$ the following
relation is valid both for ordinary Hopf algebras as well as for super-Hopf
ones:
\begin{equation}
\label{univf}
\begin{array}{c}
(m\otimes
\mbox{\rm  id})\left[ (1\otimes {\cal R}_1\otimes {\cal R}_2)(\tau \otimes
\mbox{\rm id})(\mbox{\rm  id}\otimes \tau )(\mbox{\rm  id}\otimes
\mbox{\rm  id}\otimes {\bf S}^{-1})(\Delta \otimes \mbox{\rm  id})\Delta
(e_s)\right] = \\ =(1\otimes e_s){\cal R.}
\end{array}
\end{equation}

Let us rewrite the third defining relation,
${\cal R} \Delta (e) = \tau \Delta (e) {\cal R}$, in terms of
structure constants,
\begin{equation}
\label{mult}(-1)^{\sigma _k\sigma _l+\sigma _k\sigma _j}\Delta
_i^{kl}m_{lj}^te_ke^j=(-1)^{\sigma _p\sigma _q}\Delta _i^{pl}m_{qp}^te^qe_l.
\end{equation}
Here $\sigma _k\equiv \sigma (k)$ is the grading function. From the formulas
(\ref{univf}) and (\ref{mult}) the explicit form of multiplication rules
 follows:
\begin{equation}
\label{eset}e_se^t=\sum_{n,l,k,u,j}(-1)^{\sigma _n(\sigma _l+\sigma
_k)+\sigma _u\sigma _k+\sigma _s\sigma _t}m_{nuk}^t\mu _s^{klj}({\bf S}
^{-1})_j^ne^ue_l.
\end{equation}

Despite the transparency of these rules it is not easy to use them directly.
In close analogy with the case of the ordinary double some
additional restructuring of the formula (\ref{eset}) is necessary.
Calculate two similar expressions: one for the element $e^t$,
\begin{equation}
\label{et}\Phi (e^t)\equiv
(-1)^{\sigma _u\sigma _k}m_{nuk}^te^n\otimes
e^k\otimes e^u,
\end{equation}
the other for $e_s$,
\begin{equation}
\label{es}
\begin{array}{c}
\Psi (e_s)\equiv (\tau \otimes
\mbox{\rm  id})(\mbox{\rm  id}\otimes \tau )(\mbox{\rm  id}\otimes
\mbox{\rm  id}\otimes {\bf S}^{-1})\Box (e_s)= \\ (-1)^{\sigma _l\sigma
_j+\sigma _k\sigma _j}\Box _s^{klj}({\bf S}^{-1})_j^ne_n\otimes e_k\otimes
e_l,
\end{array}
\end{equation}
with $\Box \equiv \mu (\mu \otimes $id$)$, $\mu$ -- the multiplication
in the dual Lie superalgebra ($BRST_q$ in our case). To write down
the product $e_s\cdot e^t$ it is sufficient to contract the first and
the second tensor factors and to multiply the third ones:
\begin{equation}
\label{neset}\left( -1\right) ^{\sigma_s \sigma_t}
e_s \cdot e^t=\left\langle \Phi ^{^{\prime }}(e^t),\Phi ^{^{\prime
}}(e_s)\right\rangle \left\langle \Phi ^{^{\prime \prime }}(e^t),\Phi
^{^{\prime \prime }}(e_s)\right\rangle \Phi ^{^{\prime \prime \prime
}}(e^t)\cdot \Phi ^{^{\prime \prime \prime }}(e_s).
\end{equation}

Applying these formulas to the pair $(PTSA_q,BRST_q)$ we obtain the
Hopf superalgebra $SD(PTSA_q,BRST_q)$
with the defining relations:

\begin{equation}
\label{dub1}
\begin{array}{l}
\,\left[ T,S\right] =0; \\
\,\left[ \tau ,\xi \right] =\frac h2\xi ; \\
\,\left[ S,\tau \right] =hs-2 \frac{h \xi}{\sinh \left( h \right)}
\cosh (\frac 12hT); \\
\,\left[ T,\tau \right] =0; \\
\,\left[ T,\xi \right] =0;
\end{array}
\,
\begin{array}{l}
\,\left\{ S,S\right\} =2
\frac{\sinh \left( hT\right) }{\sinh \left( h\right) }; \\ \,\left\{ \xi
,\xi \right\} =0; \\
\,\left\{ s,\xi \right\} =2\sinh \left( \frac 12hT\right) ;
\end{array}
\end{equation}

\begin{equation}
\label{dub9}
\begin{array}{l}
\Delta T=T\otimes 1+1\otimes T; \\
\Delta \xi =\xi \otimes 1+1\otimes \xi ; \\
\Delta S=e^{
\frac{hT}2}\otimes S+S\otimes e^{-\frac{hT}2}; \\ \Delta \tau =\tau \otimes
1+1\otimes \tau +\frac h{\sinh \left( h\right) }\xi \otimes \xi ;
\end{array}
\end{equation}
\begin{equation}
\label{dub13}
\begin{array}{cc}
{\bf S}(T)=-T; & {\bf S}(\tau )=-\tau ; \\ {\bf S}(S)=-S; & {\bf S}(\xi
)=-\xi .
\end{array}
\end{equation}
It is easy to check that
the universal ${\cal R}$-matrix (\ref{rmat}) realize the triangularity of
this quantum superdouble.

\section{Deformations of super Lie-Poisson structures induced by superdouble}

Applying quantum duality to the algebra $PTSA_q$ one can introduce the
canonical parameter $p$ dual to $h$ \cite{lyakhczec}. The composition
$$ \left\{
s,s\right\} =2p\frac{\sinh \left( hT\right) }{\sinh \left( h\right) }
$$
is the only relation that changes. In the $(BRST)_q$-algebra the co-product
$ \Delta (\tau )$ also aquires the dual parameter: $$ \Delta \tau =\tau
\otimes 1+1\otimes \tau +\frac{hp}{\sinh \left( h\right) } \xi \otimes \xi
; $$ (compare with (\ref{brst-def})). As a result we obtain the
two-parametric family \\ $SD^{hp}(PTSA,BRST)$ of quantum doubles.  It can
be observed that in the Hopf algebra (\ref{dub1},\ref{dub9},\ref{dub13})
the composition $ [\tau, \xi]$ allows the rescaling
$$
 [\tau, \xi] = \frac{1}{2} \alpha h \xi $$ with the additional arbitrary
parameter $\alpha$. We shall consider the case $\alpha =2$ (in order to
have the necessary classical limts) and chose the one-dimensional
family of Hoph algebras putting $p=1-h$. The defining relations for
$SD_{\alpha = 2}^{h,1-h}\equiv SD^{(h)}$ are
\begin{equation} \label{line}
\begin{array}{l} \,\left[ \tau ,\xi \right] = h\xi ; \\ \,\left\{
S,S\right\} =2\left( 1-h\right) \frac{\sinh (hT)}{\sinh (h)}; \\ \,\left\{
S,\xi \right\} =2\sinh \left( \frac{hT}2\right) ; \\ \,\left[ S,\tau
\right] =hS- \frac{2h\left( 1-h\right) }{\sinh \left( h\right) }\xi \cosh
\left( \frac{hT} 2\right) ; \\ \Delta (\tau )=\tau \otimes 1+1\otimes \tau
+ \frac{h(1-h)}{\sinh \left( h\right) }\xi \otimes \xi ; \\ \Delta \left(
S\right) =\exp \left( \frac 12hT\right) \otimes S+S\otimes \exp \left(
-\frac 12hT\right) ;
\end{array} \end{equation}
(from here on we expose only nonzero supercommutators and
nonprimitive coproducts).

According to the general theory of quantum double \cite{Sem} the elements of
the set $SD^{(h)}$ can be presented as the deformation quantizations, the
corresponding Lie superbialgebra can be constructed using the classical
Manin triple. Now we shall show that the set $SD^{(h)}$ induces deformations
of super Lie-Poisson (SL-P) structures thus attributed to the Hopf algebras
in $
SD^{(h)}$.

Consider the Hopf algebra $H^{(0)}\in SD^{(h)}$ described by the relations
(\ref{line}) in the limit $h \rightarrow 0$:
\begin{equation}
\label{0point1}
\begin{array}{l}
\,\left[ S,\tau \right] =-2\xi ; \\
\,\left\{ S,S\right\} =2T;
\end{array}
\end{equation}
\begin{equation}
\label{0point2}\Delta (\tau )=\tau \otimes 1+1\otimes \tau +\xi \otimes \xi
.
\end{equation}
This limit can be interpreted as a quantized semidirect product
$(PTSA \vdash {\rm Ab})_q$.
The corresponding analytical
\cite{lyakh-prep} variety ${\cal D}_{\mu \theta }^{(0)}$ of Hopf algebras is
defined by the compositions
\begin{equation}
\label{0facet}
\begin{array}{l}
\,\left[ S,\tau \right] =-2\mu \xi ; \\
\,\left\{ S,S\right\} =2\mu T; \\
\Delta (\tau )=\tau \otimes 1+1\otimes \tau +\theta \,\xi \otimes \xi .
\end{array}
\end{equation}
These relations correspond to the quantized SL-P structure in which
the cocommutative superalgebra $(PTSA \vdash {\rm Ab})$ is deformed in the
direction
of the Poisson bracket $ \{ \xi, \xi \} = \tau \theta$. This
quantization looks trivial, the multiplications in (\ref{0facet}) do not
depend on $\theta$.

In the opposite limit $ h \rightarrow 1$ the Hopf algebra
$H^{(1)}\in SD^{(h)}$ presents a nontrivial deformation of a
semidirect product $(BRST \vdash {\rm Ab})$:
\begin{equation}
\label{1point1}
\begin{array}{l}
\,\left[ \tau ,\xi \right] =\xi ; \\
\,\left[ S,\tau \right] =+S; \\
\,\left\{ S,\xi \right\} =2\sinh \left( \frac T2\right) ;
\end{array}
\end{equation}
\begin{equation}
\label{1point2}\Delta \left( S\right) =\exp \left( \frac 12T\right) \otimes
S+S\otimes \exp \left( -\frac 12T\right) .
\end{equation}
The procedure analogous to that used for $H^{(0)}$ leads to the analytical
variety ${\cal D}_{\mu \theta
}^{(1)}$ of Hopf algebras
\begin{equation}
\label{1facet}
\begin{array}{l}
\,\left[ S,\tau \right] =+\mu S; \\
\,\left[ \tau ,\xi \right] =\mu \xi ; \\
\,\left\{ S,\xi \right\} =2\frac \mu \theta \sinh \left(
\frac{\theta T}2\right) ; \\ \Delta \left( S\right) =\exp \left( \frac
12\theta T\right) \otimes S+S\otimes \exp \left( -\frac 12\theta T\right) .
\end{array}
\end{equation}
They have dual classical limits.
The two varieties ${\cal D}_{\mu \theta }^{(0)}$ and ${\cal D}_{\mu \theta
}^{(1)}$ intersect in the trivial point -- the Abelian and coAbelian Hopf
algebra $H_{00}^{(0)}=H_{00}^{(1)}$.

Let us show that there exists the contineous deformation \cite{lyakh-prep}
of the SP-L structure ${\cal D}_{\mu \theta }^{(0)} $ in the direction of
${\cal D} _{\mu \theta }^{(1)}$. The first order deforming functions for
such a deformation is a field on ${\cal D}_{\mu \theta }^{(0)}$ tangent to
the flow connecting ${\cal D}_{\mu \theta }^{(0)}$ and ${\cal D}_{\mu
\theta }^{(1)}$ . Evaluating the difference between the compositions
(\ref{1facet}) and (\ref{0facet}) and comparing it with the curve
(\ref{line}) as a representative of the flow we get the deforming field
${\cal F}_{\mu \theta }^{(0)}$: \begin{equation} \label{def-field}
\begin{array}{l} \,\left[ S,\tau \right] =+\mu S+2\mu \xi ; \\ \,\left[
\tau ,\xi \right] =\mu \xi ; \\ \,\left\{ S,S\right\} =-2\mu T; \\
\,\left\{ S,\xi \right\} =\mu T; \\ \Delta \left( S\right) =\frac 12\theta
\,T\wedge S; \\ \Delta (\tau )=-\theta \,\xi \otimes \xi . \end{array}
\end{equation} One can integrate the equations $$ \frac{\partial H_{\mu
,\theta }^{(h)}}{\partial h}_{\mid h=0}={\cal F}_{\mu \theta }^{(0)} $$
imposing the boundary conditions $H_{\mu ,\theta }^{(0)}\in {\cal D}_{\mu
\theta }^{(0)},\;H_{\mu ,\theta }^{(1)}\in {\cal D}_{\mu \theta
}^{(1)},\;$ and $H_{1,1}^{(h)}=SD^{(h)}$. One of the possible solutions is
the 3-dimensional variety ${\cal D}_{\mu \theta }^{(h)}$ of Hopf algebras
with compositions \begin{equation} \label{variety} \begin{array}{l}
\,\left[ S,\tau \right] =+\mu hS-\mu \frac{2h(1-h)}{\sinh \left( h\right)
}\xi \cosh \left( \frac 12h\theta T\right) ; \\ \,\left[ \tau ,\xi \right]
=\mu h\xi ; \\ \,\left\{ S,S\right\} =2\frac \mu \theta \left( 1-h\right)
\frac{\sinh (h\theta T)}{\sinh (h)}; \\ \,\left\{ S,\xi \right\} =2\frac
\mu \theta \sinh \left( \frac 12h\theta T\right) ; \\ \Delta \left(
S\right) =\exp \left( \frac 12h\theta T\right) \otimes S+S\otimes \exp
\left( -\frac 12h\theta T\right) ; \\ \Delta (\tau )=\tau \otimes
1+1\otimes \tau +\frac{h(1-h)}{\sinh \left( h\right) }\theta \,\xi \otimes
\xi . \end{array} \end{equation} For each $h^{\prime }\in \left[
0,1\right] $ fixed the 2-dimensional subvariety ${\cal D}_{\mu \theta
}^{(h^{\prime })}$ defines the SL-P structure: \begin{equation}
\label{l-p1} \begin{array}{l} \,\left[ S,\tau \right] =+\mu h^{\prime
}S-\mu \frac{2h^{\prime }(1-h^{\prime })}{\sinh \left( h^{\prime }\right)
}\xi ; \\ \,\left[ \tau ,\xi \right] =\mu h^{\prime }\xi ; \\ \,\left\{
S,S\right\} =2\mu \left( 1-h^{\prime }\right) \frac{h^{\prime }}{\sinh
(h^{\prime })}T; \\ \,\left\{ S,\xi \right\} =\mu h^{\prime }T;
\end{array} \end{equation} \begin{equation} \label{l-p2} \begin{array}{l}
\delta \left( S\right) =\frac 12h^{\prime }\theta T\wedge S; \\ \delta
(\tau )=\frac{h^{\prime }(1-h^{\prime })}{\sinh \left( h^{\prime }\right)
}\theta \,\xi \otimes \xi ; \end{array} \end{equation} described here as a
pair of superalgebra (\ref{l-p1}) and supercoalgebra (\ref{l-p2}). For
$h^{\prime }\in \left( 0,1\right) $ these structures are equivalent. But
this is not true for the limit points -- ${\cal D}_{\mu \theta }^{(0)}$
and ${\cal D}_{\mu \theta }^{(1)}$ represent two different contractions of
the quantized SL-P structure ${\cal D}_{\mu \theta }^{(h^{\prime })}\mid
_{h^{\prime }\in \left( 0,1\right) }$. Thus the main statement is proved:
the SL-P structure (\ref{0point1},\ref{0point2}) (``trivially'' quantized
as ${\cal D}_{\mu \theta }^{(0)}$) can be deformed in the direction of
Hopf algebras belonging to ${\cal D}_{\mu \theta }^{(1)}$ (that is -- by
the field ${\cal F}_{\mu \theta }^{(0)}$ ) to obtain the quantization
\begin{equation} \label{newquant} \begin{array}{l} \,\left[ S,\tau \right]
=+\mu hS-\mu \frac{2h(1-h)}{\sinh \left( h\right) }\xi \cosh \left( \frac
12 h^2 T\right) ; \\ \,\left[ \tau ,\xi \right] =\mu h\xi ; \\ \,\left\{
S,S\right\} =2\frac \mu h \left( 1-h\right) \frac{\sinh (h^2 T)}{\sinh
(h)}; \\ \,\left\{ S,\xi \right\} =2\frac \mu h \sinh \left( \frac 12h^2
T\right) ; \\ \Delta \left( S\right) =\exp \left( \frac 12h^2 T\right)
\otimes S+S\otimes \exp \left( -\frac 12h^2 T\right) ; \\ \Delta (\tau
)=\tau \otimes 1+1\otimes \tau +\frac{h^2 (1-h)}{\sinh \left( h\right) }
\,\xi \otimes \xi . \end{array} \end{equation} One of the classical limits
(for $\mu \rightarrow 0$) lay in the facet
 ${\cal D}_{0 \theta }^{(h)} $ of classical supergroups (\ref{dub9}).
Note that depite these properties the Hopf algebra (\ref{newquant})
is a quantization of the same super Lie bialgebra as in the trivial
canonical quantization of the proper time group cotangent bundle
(\ref{0facet}). This is easily checked by evaluating the first order
terms in the expansion of the compositions (\ref{newquant}) with respect to
$\mu$ and $h$.
This deformation is induced by the quantum superdouble construction.

Earlier (see \cite{lyakh-prep}) it was demonstrated that
quantum double could induce even more complicated deformations of L-P
structures where the corresponding groups and algebras of observables are
not only deformed but also quantized. In the case disscussed above the
procedure presented in \cite{lyakh-prep} does not lead to nontrivial results
The variety ${\cal
D}_{\mu \theta }^{(0)}$ lifted in the domain of non(anti)commutative and
nonco(anti)commutative Hopf algebras will have edges equivalent to its
internal points. This is a consequence of the equivalence of all the
Hopf algebras corresponding to the internal points of
${\cal D}_{\mu \theta }^{(h)}$.

\section{Conclusions}
Analyticity plays the important role in the selection of admissable
transformations of Poisson structures.
Although the SL-P structures corresponding to
$ \{ {\cal D}_{\mu \theta }^{(h^{\prime })} \mid  h' \in (0,1) \} $
are equivalent,
the contineous ''rotation'' of
${\cal D}_{\mu \theta }^{(h^{\prime })}$ breaks the analiticity.
This is in accordance with the fact that the
compositions (\ref{l-p1}, \ref{l-p2}) with different $ h^{\prime}$'s
do not form super Lie bialgebra. This
effect was first observed in \cite{lyakh-pap} for a nonsuper case.

The deformations ${\cal D}_{\mu \theta }^{(0)}\longrightarrow
{\cal D}_{\mu \theta }^{(h^{\prime })}$ might be of considerable physical
importance. We would like to stress that in these deformations both the
supergroup and the Poisson superalgebra of its coordinate functions are
deformed simultaneously. Moreover, the process can not be subdevided into
successive deformations of group and algebra for the reasons
described above. Thus the deformation
of the dynamics must be accompanied by the deformation of the geometry.
In our particular case the Lie superalgebra of the cotangent bundle
$T^*(PTSG)$ can be quantized (retaining the Hopf structure)
if the Abelian subalgebra of the cotangent space is simultaneously
deformed into the $BRST$-like algebra and one of
the canonical classical limits becomes isomorphic to the classical
double of $PTSG$ and $BRST$ groups.

It must be mentioned that other
methods of unification such as crossproducts or cocyclic cross- and
bicrossproducts of Hopf algebras do not lead to nontrivial algebraic
constructions in the case of $PTSA_q$ and $BRST_q$.

{\large \bf Aknowledgments}

One of the authors (V.D.L.) would like to express his gratitude to
colleagues in the Institute of Physics of the University of Guanajuato for
their warm hospitality during the completion of this work.

Supported in part by the Russian Foundation for Fundamental Research
 grant N 97-01-01152 (V.D.L.) and CONACyT grant 3898P-E9608 (V.I.T.).

\end{document}